\documentclass[A4]{article} 
\usepackage{color,mathrsfs}
\usepackage{amssymb,amsmath,latexsym}
\usepackage{graphicx,amssymb}
\usepackage[utf8]{inputenc}
\usepackage{authblk}
\usepackage[T1]{fontenc}

\newtheorem{thm}{Theorem}
\newtheorem{prop}{Proposition}

\newtheorem{defin}{Definition}

\newcommand{\N}{\mathbb N}
\newcommand{\R}{\mathbb R}
\newcommand{\Z}{\mathbb Z}

\newcommand{\erre}{\mathbb R}

\def\tilde{\widetilde}
\def\epsilon{\varepsilon}

\def\boxit#1{\vbox{\hrule\hbox{\vrule\kern.75truemm
\vbox{\kern.75truemm#1\kern1truemm}\kern1truemm\vrule}\hrule}}

\def\qed{\hskip 1mm\boxit{}\hskip 1mm}

\title{Convergence in variation  for the \\ multidimensional generalized sampling series and applications to smoothing\\ for digital image processing}

\author{ {\bf Laura Angeloni, \hskip0.3cm Danilo Costarelli, \hskip0.3cm Gianluca Vinti} \\
Department of Mathematics and Computer Science\\
University of Perugia, Via Vanvitelli 1, 06123, Perugia, Italy\\
{\small {\tt laura.angeloni@unipg.it} \hskip0.3cm - \hskip0.3cm {\tt danilo.costarelli@unipg.it}}\\
{\small {\tt gianluca.vinti@unipg.it}}}

\date{}

\begin{document}
\maketitle

\begin{abstract}
In this paper we study the problem of the convergence in variation for the generalized sampling series based upon averaged-type kernels in the multidimensional setting. As a crucial tool, we introduce a family of operators of sampling-Kantorovich type for which we prove convergence in $L^p$ on a subspace of $L^p(\R^N)$: therefore we obtain the convergence in variation for the multidimensional generalized sampling series by means of a relation between the partial derivatives of such operators acting on an absolutely continuous function $f$ and the sampling-Kantorovich type operators acting on the partial derivatives of $f$. Applications to digital image processing are also furnished.
\end{abstract}

\section{Introduction} \label{sec-intro} 

In this paper we present approximation results in BV-spaces for the generalized sampling series in the multidimensional frame.

   The above sampling series, defined as
$$
\mbox{(I)} \hskip1cm (S_w f)(t):=\sum_{k\in\Z} f\left({k\over w}\right)\chi(w{ t}-{k}), \ \ { t}\in\R,\ w>0,
$$
has been introduced by P.L. Butzer (see, e.g., \cite{BUSPST}) and revealed to be very interesting for both the theoretical and applicative aspects. Indeed, several approximation results have been obtained in the last forty years with respect to different notions of convergence, such as uniform, $L^p$, and modular convergence (\cite{BV97,MV03,BABUSTVI2}). Moreover, they have many important applications to Signal Theory since they provide an approximate sampling formula, which allows to reconstruct not necessarily band-limited signals. By its multivariate generalization, that is,
$$
\mbox{(II)} \hskip1cm (S_w f)({\tt t}):=\sum_{{\tt k}\in\Z^N} f\left({{\tt k}\over w}\right)\chi(w{\tt t}-{\tt k}), \ \ {\tt t}\in\R^N,\ w>0,
$$
several problems of image processing can be treated (see, e.g., \cite{BMSVV}). 

A natural setting to study some issues involving digital images is furnished by the spaces of functions of bounded variation (\cite{GO1,SC1,VA1}). In this direction, it is interesting to have results about estimates and convergence in variation for the above discrete operators. In the case of functions of one-variable, the problem has been faced in \cite{ACV1} for the operators (I) based upon averaged type kernels.

    Here, also with the aim to consider applicative aspects, we deal with the multivariate case, namely we study results about approximation in variation by means of the operators (II) based upon multidimensional product kernels of averaged type, namely
$$
\bar \chi_m({\tt t}):=\prod_{i=1}^N \bar\chi_{i,m}(t_i)
$$
where 
$$
\bar\chi_{i,m}(t):={1\over m}  \int_{-{m\over 2}}^{m\over 2} \chi_i(t+v)\,dv,
$$
for some $m\in\N$, and $\chi_i:\R\longrightarrow\R$ are suitable one-dimensional kernels. We will use the multidimensional generalization of the classical Jordan variation introduced by Tonelli (see, e.g., \cite{TO}) and later generalized by T. Rad\'o (\cite{RA}) and C. Vinti (\cite{VI}) to the case of functions of $N$-variables.
    
    In particular, we establish a variation diminishing-type property (Proposition \ref{var-dim}) and a convergence theorem (Theorem \ref{conv_var}).
    
    Due to form of the above discrete sampling type operators, to obtain results about approximation in variation is a very delicate problem. The strategy we propose is to introduce a new family of Kantorovich type operators defined as 
$$
\mbox{(III)} \hskip0.4cm (K_{w,j} f)({\tt t}):=\sum_{{\tt k}\in\Z^N} \left[w\int_{k_j\over w}^{k_{j}+1\over w}f\left({k_1\over w},\ldots, u\ldots, {k_N\over w}\right)\,du\ \right] \chi(w{\tt t}-{\tt k}), 
$$
${\tt t}\in\R^N$, $w>0$, $j=1,\ldots,N$, for which we establish $L^p$-convergence in a subspace of $L^p(\R^N)$ (Theorem \ref{conv_p}), therefore obtaining a convergence result for the new class of operators (III). We notice that the above result is obtained through the use of the $\tau-$modulus of continuity (see \cite{K69,DS76,SP}) which seems to be the most suitable approach in this setting. 
    
    The crucial point in order to reach the convergence in variation is to find a link between the two classes of operators: to this aim we prove (Proposition \ref{prop-der}) a relation between the gradient of the multivariate generalized sampling series of a function $f$
and the family of Kantorovich-type discrete operators acting on the partial derivatives of $f$.

Some of the results of the present paper may have an applicative interpretation: in particular, the variation diminishing type property (Proposition \ref{var-dim}), ensuring that the generalized sampling series (II) have smaller variation than the function on which they act, can be viewed as a smoothing procedure. Indeed, if the function $f$ is a gray scale digital image, the action of the operators and the consequent variation diminishing corresponds to reducing the jumps of gray levels with respect to the original image, hence producing a smoothing effect. In Section \ref{sec6} we discuss such digital image processing applications in details, furnishing some numerical examples.


\section{Notations and preliminaries}\label{sec-not}

Our results will be set in the frame of the space of multivariate functions of bounded variation: in particular we will consider the concept of variation introduced by Tonelli (\cite{TO}) for two variables, extended to the general case of $\R^N$ by Rad\'o and C. Vinti (\cite{RA,VI}). Here we will  recall it.  

For a function $f:\R^N \rightarrow \R$ and ${\tt
x}=(x_1,\dots,x_N)\in\R^N$, we will use the notation
$$
{\tt x}'_j= (x_1,\dots, x_{j-1}, x_{j+1},\dots,x_N)\in\erre^{N-1},\ \
{\tt x}=({\tt x}'_j,x_j),\ \ f({\tt x})=f({\tt x}'_j,x_j),
$$
if we are interested in
the $j-$th coordinate of ${\tt x}$,  $j=1,\dots ,N$. Moreover, given an $N-$dimensional interval $I=\prod_{i=1}^N
[a_i,b_i],$ by $I'_j=[{\tt a}'_j,{\tt b}'_j]$ we will denote the
\mbox{$(N-1)-$dimensional} interval obtained deleting by $I$ the
$j-$th coordinate, i.e.,
$$
I=[{\tt a}'_j,{\tt b}'_j]\times[a_j,b_j],\ \ j=1,\dots ,N.
$$ 
Given a vector ${\tt x}\in\R^N$ and $\alpha \in \R$, we will use the usual notation for products and quotients, i.e., ${\alpha{\tt x}}=\left({\alpha x_1 },\ldots,{ \alpha x_N}\right)$ and, for $\alpha \neq 0$,
${{\tt x}\over \alpha}=\left({x_1\over \alpha},\ldots,{x_N\over \alpha}\right)$.

By $R_{loc}(\R^N)$ we will denote the space of locally Riemann integrable functions on $\R^N$, while $M(\R^N)$ will denote the space of all the measurable and bounded functions $f:\R^N\rightarrow \R$.

\begin{defin}
A function $f \in M(\R^N)$ is said to be {\it of bounded variation} ($f\in BV(\R^N)$) if $V_{\R}[f({\tt x}'_j,\cdot)]$
(the usual Jordan one-dimensional variation of the $j-$th section of $f$) is finite for a.e. ${\tt x}'_j \in\R^{N-1}$ and
$$
\int_{\R^{N-1}} V_{\R}[f({\tt x}'_j,\cdot)] \,d{\tt x}'_j <+\infty,
$$ 
for every $j=1,\ldots,N$.
\end{defin}
For more details about $BV$-spaces, see, e.g., \cite{AVI-2007,A-2013,AV-2013,AV-2014,AD-2014,AP1}.

Let us now recall how to compute the multidimensional Tonelli variation. The first step is to consider,
for $I=\prod_{i=1}^N
[a_i,b_i]$ and $j=1,\dots,N$ the $(N-1)-$dimensional integrals
$$
\Phi_j(f,I):= \int_{[{\tt a}'_j,{\tt b}'_j]} V_{[a_j,b_j]}[f({\tt x}'_j,\cdot)] 
d{\tt x}'_j,
$$
where $V_{[a_j,b_j]}[f({\tt x}'_j,\cdot)]$ is the usual one-dimensional (Jordan)
variation of the $j-th$ section of $f$. 
Let now $\Phi(f,I)$ be the Euclidean norm of the vector
$(\Phi_1(f,I),\dots, \Phi_N(f,I))$, namely
$$
\Phi(f,I):= \left\{\sum_{j=1}^N \Phi_j^2(f,I)\right\}^{1\over 2},
$$
where $ \Phi(f,I)=+\infty$ if $\Phi_j(f,I)=+\infty$ for some
$j=1,\dots,N$.

Then the variation of $f$ on $I\subset \R^N$ is defined 
as
$$
V_I[f] := \sup \sum_{k=1}^m \Phi(f,J_k),
$$
where the supremum is taken over all the finite families of
$N-$dimensional intervals $\{J_1,\dots,J_m\}$ which form
partitions of $I$.

Passing to the supremum over all the intervals $I\subset
\R^N$, we obtain the variation of $f$ over the whole $\R^N$, i.e.,
$$
V[f]:= \sup_{I\subset \R^N} V_I[f].
$$
It is well known that, for every $f\in BV(\erre^N)$,
$\nabla f$ exists a.e. in $\R^N$ and ${\partial f \over \partial x_j} \in
L^1(\R^N)$, for every $j=1,\ldots,N$ (see e.g. \cite{RA,VI}).

We now recall the notion of absolute continuity in the sense of Tonelli. 

\begin{defin}
A function $f:\R^N\rightarrow \R$ is locally absolutely continuous in the sense of Tonelli ($f\in AC_{loc}(\R^N)$) if, for every interval
$I=\prod_{i=1}^N [a_i,b_i]$ and for every $j=1,2,\dots,N$, the $j-$th
section of $f$, $f({\tt x}'_j,\cdot):[a_j,b_j]\rightarrow \R$ is absolutely continuous for
almost every ${\tt x}'_j\in [{\tt a}'_j,{\tt b}'_j]$.
\end{defin} 

It is a well known result that, if $f\in
BV(\R^N)\cap AC_{loc}(\R^N),$ then
$$
V[f]=\int_{\R^N} |\nabla f({\tt x})| \,d{\tt x}
$$
(see \cite{RA,VI,GIU,AV-2015}), that is, an integral representation for the variation of $f$ holds.

We will therefore denote by $AC(\R^N)$ the space of all the functions
$f\in BV(\R^N) \cap AC_{loc}(\R^N)$ ({\it absolutely continuous functions}).

We will study approximation results for the {\it multivariate generalized sampling series}, namely a family of discrete operators defined as
$$
(S_w f)({\tt t}):=\sum_{{\tt k}\in\Z^N} f\left({{\tt k}\over w}\right)\chi(w{\tt t}-{\tt k}), \ \ {\tt t}\in\R^N,\ w>0:
$$
such operators are the multidimensional version of the generalized sampling series (see, e.g., \cite{BFS1,BFS2,BMSVV}).

Here $\chi$ is a kernel, that is, a function $\chi:\R^N\rightarrow \R$ that satisfies the following assumptions:

\begin{description} 
\item{$(\chi_1)$}
$\chi \in L^1(\R^N)$ is 
such that $\sum_{{\tt k}\in\Z^N} \chi({\tt u-k})=1$, for every ${\tt u}\in\R^N$;
\item{$(\chi_2)$} 
$A_{\chi}:=\sup_{{\tt u}\in\R^N}\sum_{{\tt k}\in\Z^N} |\chi({\tt u-k}) |<+\infty$, where the convergence of the series is uniform on the compact sets of $\R^N$. 
\end{description}

The above assumptions are quite standard working with discrete families of operators: see e.g., \cite{BFS1}.
We point out that the operators $(S_w f)_{w>0}$ are well-defined, for example, for every $f\in BV(\R^N)$: indeed in this case $f$ is bounded and since $|f({\tt t})| \le M$, for some $M>0$ and for every ${\tt t}\in\R^N$, we have
$$
|(S_w f)({\tt t})|\le M\sum_{{\tt k}\in\Z^N} \chi(w{\tt t}-{\tt k}) \le MA_{\chi},
$$
by $(\chi_2)$.

In particular, in the present paper  we will study the convergence in variation for the multivariate generalized sampling series with product kernels of {\it averaged type}, that is kernels of the form 
\begin{equation}\label{av_ker}
\bar\chi_m({\tt t}):=\prod_{i=1}^N \bar\chi_{i,m}(t_i)
\end{equation}
where 
\begin{equation*}
\bar\chi_{i,m}(t):={1\over m}  \int_{-{m\over 2}}^{m\over 2} \chi_i(t+v)\,dv,
\end{equation*}
for some $m\in\N$, and $\chi_i:\R\longrightarrow\R$ is a (one-dimensional) kernel for every $i=1,\ldots,N$ (i.e., satisfying $(\chi_1)$ and $(\chi_2)$ with $N=1$). 

Notice that it is easy to see that $\bar\chi_m$ is a kernel itself and 
moreover, by the Fubini-Tonelli theorem, 
\begin{equation}\label{norm_kernel}
\begin{split}
\Vert \bar\chi_{i,m} \Vert_{L^1(\R)}& = \int_{\R} \left| {1\over m}  \int_{-{m\over 2}}^{m\over 2} \chi_i(t+v)\,dv \right|\,dt 
\\ & \le  {1\over m}  \int_{-{m\over 2}}^{m\over 2} \int_{\R} | \chi_i(t+v)| \,dt\,dv = \Vert \chi_i\Vert_{L^1(\R)}
\end{split}
\end{equation}
and, by this, 
\begin{equation*}
\Vert \bar\chi_{m} \Vert_1  \le \prod_{i=1}^N \Vert \chi_i\Vert_{L^1(\R)}.
\end{equation*}

We point out that is it easy to give examples of product kernels of averaged type (see Section \ref{sec6}) to which our results can be applied.

The corresponding multivariate generalized sampling series (associated to the averaged product kernel $\bar \chi_m$) will be denoted as
$$
(\bar S^m_w f)({\tt t}):=\sum_{{\tt k}\in\Z^N} f\left({{\tt k}\over w}\right)\bar\chi_m(w{\tt t}-{\tt k}), \ \ {\tt t}\in\R^N,\ w>0.
$$

It is easy to see that $\bar\chi_m$ is differentiable and, taking into account that, obviously, 
\begin{equation*}\label{av_ker_der}
{\partial\bar\chi_m \over \partial t_j} ({\tt t}) ={1\over m} \prod_{i\neq j} \bar\chi_{i,m}(t_i)\left[\chi_j\left(t_j+{m\over 2}\right)-\chi_j\left(t_j-{m\over 2}\right)\right], \ \ {\tt t}\in\R^N,
\end{equation*}
it is possible to write, for every $j=1,\ldots,N$, ${\tt t}\in\R^N$, $w>0$, 
\begin{equation}\label{sampl_der}
{\partial \bar S^m_w f \over \partial t_j} ({\tt t}) = {w\over m} \sum_{{\tt k}\in\Z^N} f\left({{\tt k}\over w}\right)\prod_{i\neq j} \bar \chi_{i,m}(wt_i-k_i)\left[\chi_j\left(wt_j-k_j+{m\over 2}\right)-\chi_j\left(wt_j-k_j-{m\over 2}\right)\right].
\end{equation}
Notice that, since $f$ is bounded, ${\partial \bar S^m_w f \over \partial t_j} ({\tt t})$ exists for every ${\tt t}\in\R^N$ since 
\begin{equation}\label{sampl_der2}
\begin{split}
\left|{\partial \bar S^m_w f \over \partial t_j} ({\tt t})\right| &\le {w\over m} M \sum_{{\tt k}\in\Z^N} \prod_{i\neq j} |\bar \chi_{i,m}(wt_i-k_i)| \left[\left|\chi_j\left(wt_j-k_j+{m\over 2}\right)\right|+\left|\chi_j\left(wt_j-k_j-{m\over 2}\right)\right|\right] \\
& \le {2w\over m} M  \prod_{i=1}^NA_{\chi_i},
\end{split}
\end{equation}
again by $(\chi_2)$ for each one-dimensional kernel $\chi_i$.

One of the main ideas of this paper is to establish a relation between the gradient of the multivariate generalized sampling series of a function $f$
and a family of Kantorovich-type discrete operators, that we now introduce, acting on the partial derivatives of $f$. This generalizes the analogous result, in the one-dimensional case, obtained in \cite{ACV1}, similarly to what happens between the Bernstein polynomials and their Kantorovich version (see, e.g., \cite{A-2001}).  

We therefore introduce the following family of multidimensional discrete operators of sampling-Kantorovich type:
$$
(K_{w,j} f)({\tt t}):=\sum_{{\tt k}\in\Z^N} \left[w\int_{k_j\over w}^{k_{j}+1\over w}f\left({k_1\over w},\ldots, u\ldots, {k_N\over w}\right)\,du\ \right] \chi(w{\tt t}-{\tt k}), \ \ {\tt t}\in\R^N,\ w>0,
$$
$j=1,\ldots,N$.
Notice that, similarly to the case of $(S_w f)_{w>0}$, it is easy to see that the operators $(K_{w,j} f)_{w>0}$, $j=1,\ldots,N$, are well-defined if, for example $f\in BV(\R^N)$, since $f$ is in particular bounded.


\section{Results for the multidimensional generalized sampling series}\label{sec-mres}

We will first prove that the operators $(\bar S^m_w f)_{w>0}$ map $BV(\R^N)$ into $AC(\R^N)$. Such result is a kind of ''asymptotic'' variation diminishing property: indeed, choosing $m\in\N$ sufficiently large, we obtain the classical variation diminishing property for the multidimensional generalized sampling series.

\begin{prop} \label{var-dim}
Let $f\in BV(\R^N)$. Then $\bar S^m_w f\in AC(\R^N)$, for every $w>0$, $m\in\N$, and 
\begin{equation} \label{variation-diminishing}
V[\bar S^m_w f]\le {1\over m} \prod_{i=1}^N\Vert \chi_i \Vert_1 V[f].
\end{equation}
\end{prop}

\noindent {\bf Proof.} Let us fix $w>0$ and $m\in\N$.  By (\ref{sampl_der}) and (\ref{sampl_der2}), the partial derivatives of $\bar S^m_w f$ exist for every ${\tt t}\in \R^N$, $j=1\ldots,N$, and are bounded: this implies that the sections of the function $f$ are locally absolutely continuous, namely $f\in AC_{loc} (\R^N)$. 
Let us fix $I=\displaystyle\prod_{i=1}^N[a_i,b_i]\subset\R^N$ and let $\{J_1,\dots,J_p\}$ be a partition of $I$,
with $J_k=\displaystyle\prod_{j=1}^N[^{(k)}a_j, ^{(k)}b_j]$, $k=1,\dots p$.
For every fixed $j=1,\dots N$, $k=1,\dots p$, similarly to the proof of Proposition 1 of \cite{ACV1}, it is possible to show that, since $ (\bar S^m_w f) ({\tt t}'_j,\cdot)\in AC_{loc}(\R)$, a.e. ${\tt t}'_j\in [^{(k)}{\tt a}'_j, ^{(k)}{\tt b}'_j]$, 
$$
V_{[^{(k)}a_j, ^{(k)}b_j]}[(\bar S^m_w f )({\tt t}'_j,\cdot)] \le {1\over m} V_{[^{(k)}a_j, ^{(k)}b_j]}[f ({\tt t}'_j,\cdot)] \prod_{i\neq j}\Vert \bar \chi_{i,m} \Vert_1 \Vert \chi_j\Vert_1
$$
and therefore, by (\ref{norm_kernel}), 
\begin{align*}
\Phi_j(\bar S^m_w f, J_k) & = \int_{[^{(k)}{\tt a}'_j, ^{(k)}{\tt b}'_j]} V_{[^{(k)}a_j, ^{(k)}b_j]}[(\bar S^m_w f) ({\tt t}'_j,\cdot)]\,d{\tt t}'_j \\ & \le {1\over m} \int_{[^{(k)}{\tt a}'_j, ^{(k)}{\tt b}'_j]} V_{[^{(k)}a_j, ^{(k)}b_j]}[f ({\tt t}'_j,\cdot)] \,d{\tt t}'_j \prod_{i=1}^N\Vert \chi_i \Vert_1 \\ & = {1\over m} \Phi_j( f, J_k) \prod_{i=1}^N\Vert \chi_i \Vert_1. 
\end{align*}

This implies that, for every $k=1,\dots,p$,
\begin{align*}
\Phi(\bar S^m_w f, J_k) &=\left\{\sum_{j=1}^N [\Phi_j(\bar S^m_w  f, J_k)]^2\right\}^{1\over 2} \le \left\{\sum_{j=1}^N 
 {1\over m^2} [\Phi_j(f, J_k)]^2 \prod_{i=1}^N\Vert \chi_i \Vert_1^2\right\}^{1\over 2} 
\\ &= {1\over m} \left\{\sum_{j=1}^N 
 [\Phi_j(f, J_k)]^2 \right\}^{1\over 2} \prod_{i=1}^N\Vert  \chi_i \Vert_1 = {1\over m} \prod_{i=1}^N\Vert \chi_i \Vert_1 \Phi(f, J_k). 
\end{align*}
Summing now over $k=1,\dots,p$ and passing to the supremum over
all the possible partitions $\{J_1,\dots,J_p\}$ of $I,$ we conclude that
$$
V_I[\bar S^m_w f]  \le {1\over m} \prod_{i=1}^N\Vert  \chi_i \Vert_1 V_I[f] 
$$
and hence, by the arbitrariness of $I\subset \R^N$,
$$
V[\bar S^m_w f]  \le {1\over m} \prod_{i=1}^N\Vert  \chi_i \Vert_1 V[f]. 
$$
This implies that $\bar S^m_w f\in BV(\R^N)$ and hence $\bar S^m_w f\in AC(\R^N)$.
\hfill\qed

We will now give the result that establishes a relation between the partial derivatives of the multidimensional sampling series and the multidimensional sampling-Kantorovich type operators acting on the partial derivatives of the function.

\begin{prop}\label{prop-der}
If $f\in AC(\R^N),$ then for every ${\tt t}\in \R^N$, $j=1,\ldots,N$,
$$
{\partial \bar S^m_w f \over \partial t_j} ({\tt t})\ =\ {1\over m} \sum_{i=1}^m \left(K_{w,j} {\partial f \over \partial t_j}\right)\left({\tt t}'_j,t_j-{m-2(i-1)\over 2w}\right),
$$
$w>0$, $m\in\N$. 
\end{prop}

\noindent {\bf Proof.}
By (\ref{sampl_der}), since $f$ is locally absolutely continuous, we have that 
\begin{align*}
{\partial \bar S^m_w f \over \partial t_j} ({\tt t}) & = {w\over m} \sum_{{\tt k}\in\Z^N} f\left({{\tt k}\over w}\right)\prod_{i\neq j} \bar\chi_{i,m}(wt_i-k_i)\left[\chi_j\left(wt_j-k_j+{m\over 2}\right)-\chi_j\left(wt_j-k_j-{m\over 2}\right)\right]
\\ &= {w \over m} \sum_{{\tt k}\in\Z^N} \left[ \int_{0}^{k_j\over w} {\partial f \over \partial t_j} \left({{\tt k}'_j \over w},u\right)\,du + f\left({{\tt k}'_j \over w},0\right)\right] \prod_{i\neq j} \bar\chi_{i,m}(wt_i-k_i)\left[\chi_j\left(wt_j-k_j+{m\over 2}\right)  \right. \\ &- \left.\chi_j\left(wt_j-k_j-{m\over 2}\right)\right] \\ 
&=  {w \over m} \sum_{{\tt k}\in\Z^N} \left[ \int_{0}^{k_j\over w} {\partial f \over \partial t_j} \left({{\tt k}'_j \over w},u\right)\,du + f\left({{\tt k}'_j \over w},0\right)\right] \prod_{i\neq j} \bar\chi_{i,m}(wt_i-k_i) \chi_j\left(wt_j-k_j+{m\over 2}\right)
 \\ & - {w \over m} \sum_{{\tt k}\in\Z^N} \left[ \int_{0}^{k_j\over w} {\partial f \over \partial t_j} \left({{\tt k}'_j \over w},u\right)\,du + f\left({{\tt k}'_j \over w},0\right)\right] \prod_{i\neq j} \bar \chi_{i,m}(wt_i-k_i) \chi_j\left(wt_j-k_j-{m\over 2}\right) .
\end{align*}
Now putting in the first series $\tilde k_j=k_j-m$ and $\tilde k_i=k_i$ for every $i\neq j$, there holds 
\begin{align*}
{\partial \bar S^m_w f \over \partial t_j} ({\tt t}) & =  
{w \over m} \sum_{\tilde {\tt k}\in\Z^N} \left[ \int_{0}^{\tilde k_j+m\over w} {\partial f \over \partial t_j} \left({\tilde{\tt k}'_j \over w},u\right)\,du + f\left({\tilde{\tt k}'_j \over w},0\right)\right] \prod_{i\neq j} \bar \chi_{i,m}(wt_i-\tilde k_i) \chi_j\left(wt_j-\tilde k_j-{m\over 2}\right)
 \\ & - {w \over m} \sum_{{\tt k}\in\Z^N} \left[ \int_{0}^{k_j\over w} {\partial f \over \partial t_j} \left({{\tt k}'_j \over w},u\right)\,du + f\left({{\tt k}'_j \over w},0\right)\right] \prod_{i\neq j} \bar \chi_{i,m}(wt_i-k_i) \chi_j\left(wt_j-k_j-{m\over 2}\right)  \\
&=  {w \over m} \sum_{{\tt k}\in\Z^N} \int_{k_j\over w}^{k_j+m\over w} {\partial f \over \partial t_j} \left({{\tt k}'_j \over w},u\right)\,du  \prod_{i\neq j} \bar \chi_{i,m}(wt_i-k_i) \chi_j\left(wt_j-k_j-{m\over 2}\right) \\
&= {w \over m} \sum_{{\tt k}\in\Z^N} \left(\sum_{i=1}^m\int_{k_j+i-1\over w}^{k_j+i\over w} \right){\partial f \over \partial t_j} \left({{\tt k}'_j \over w},u\right)\,du  \prod_{i\neq j} \bar \chi_{i,m} (wt_i-k_i) \chi_j\left(wt_j-k_j-{m\over 2}\right) \\
&= {1 \over m} \left\{ \left(K_{w,j} {\partial f \over \partial t_j} \right)\left({\tt t}'_j, t_j-{m\over 2w}\right)+ \ldots +\left( K_{w,j} {\partial f \over \partial t_j} \right)\left({\tt t}'_j, t_j-{m-2(m-1)\over 2w}\right)\right\}
\\ &= {1\over m} \sum_{i=1}^m \left(K_{w,j} {\partial f \over \partial t_j}\right)\left({\tt t}'_j,t_j-{m-2(i-1)\over 2w}\right).
\end{align*}
\hfill\qed


\section{Convergence in $L^p$ for the multidimensional sampling-Kantorovich type operators} \label{sec4}

We now study the problem of the convergence in $L^p$ for the multidimensional sampling-Kantorovich type operators that we introduced.
Such result will be also fundamental in order to prove the convergence in variation for $(\bar S_w f)_{w>0}$.

We recall that convergence in $L^p$ holds for the multidimensional generalized sampling series $(S_wf)_{w>0}$ (see \cite{BMSVV}) assuming that $f$ belongs to a suitable subspace of $L^p(\R^N)$, namely $\Lambda^p$. Since the definition of our operators $(K_{w,j} f)_{w>0}$ is very close to that one of $(S_wf)_{w>0}$, it is natural to expect that the $L^p-$convergence holds within the same subspace, that is what we will now prove. 

Let us recall the definition of $\Lambda^p$ (see \cite{BMSVV}). 

Let us consider an admissible partition over the $i-$th axis, i.e., $\Sigma_i:=(x_{i,j_i})_{j_i\in \Z}$ such that
$$
0<\underline{\Delta}:=\min_{i=1,\ldots,N} \inf_{j_i\in \Z} (x_{i,j_i}-x_{i,j_i-1}) \le \max_{i=1,\ldots,N} \sup_{j_i\in \Z} (x_{i,j_i}-x_{i,j_i-1})=: \overline{\Delta}<+\infty.
$$
We say that $\Sigma=({\tt x}_{\tt j})_{{\tt j}\in\Z^N}\subset \R^N$, ${\tt x}_{\tt j}=(x_{1,j_1},\ldots,x_{N,j_N})$, ${\tt j}=(j_1,\ldots,j_N)\in\Z^N$, is an admissible sequence if it is the cartesian product of admissible partitions $\Sigma_i=(x_{i,j_i})_{j_i\in\Z}$.

For a fixed admissible sequence $\Sigma$, the $l^p(\Sigma)-$ norm of a function $f:\R^N\rightarrow \R$ is defined as 
$$
\Vert f\Vert_{l^p(\Sigma)} := \left\{\sum_{{\tt j}\in\Z^N} \sup_{{\tt x}\in Q_{\tt j}}|f({\tt x})|^p \Delta_{\tt j}\right\}^{1\over p}  , \ \ 1\le p<+\infty
$$
where $Q_{\tt j}=\prod_{i=1}^N [x_{i,j_i-1}-x_{i,j_{i}}[$ and $\Delta_{\tt j}:=\prod_{i=1}^N(x_{i,j_i}-x_{i,j_{i}-1})$ is the volume of  $Q_{\tt j}$.
With such notations the subspace $\Lambda^p$ is defined as
$$
\Lambda^p:=\{f\in M(\R^N):\ \Vert f\Vert_{l^p(\Sigma)} <+\infty, \ \hbox{for every admissible sequence $\Sigma$}\}.
$$
 In \cite{BMSVV} it is proved that $\Lambda^p$ is a proper linear subspace of $L^p(\R^N)$, together with other properties concerning such space. In particular, we recall the following important result of convergence in $L^p$ for the $\tau-$modulus of smoothness (see \cite{K69,DS76,SP})
$$
\tau_{r}(f;\delta, M(\R^N))_p: = \Vert \omega_r(f;\cdot,\delta, M(\R^N))\Vert_p , 
$$
where 
$$
\omega_r(f;{\tt x},\delta, M(\R^N)) := \sup\left\{|\Delta^r_{\tt h} f({\tt t})|:\ {\tt t}, {\tt t+h}r \in \prod_{i=1}^N\left[x_i-{\delta r\over 2},x_i+{\delta r\over 2}\right]\right\},
$$ 
$\delta>0$, and $\Delta^r_{\tt h} f({\tt t}):=\sum_{j=0}^r(-1)^{r+j}\left(\begin{array}{l} r \\ j\end{array}\right)f({\tt t}+j{\tt h})$ denotes the differences of order $r$ at ${\tt t}\in\R^N$ with increment ${\tt h}\in\R^N$.

\begin{prop}[Prop. 7 of \cite{BMSVV}]  \label{tau_p_conv}
If $f\in \Lambda^p\cap R_{loc}(\R^N)$, $ 1\le p<+\infty$, $r\in\N$, then
$$
\lim_{\delta \to 0^+} \tau_{r}(f;\delta, M(\R^N))_p=0.
$$
\end{prop}

We point out that, of course, $\Lambda^p$ is a non trivial subspace since it contains, for example, all the functions in $M(\R^N)$ with compact support. 

\vskip0.1cm

We are now ready to state the result of convergence in $L^p$ for the multidimensional sampling-Kantorovich type operators. 

Let us now assume that $\chi$ is a kernel with compact support, i.e., $\chi({\tt t
})=0$ if ${\tt t}\not\in [-T,T]^N$, $T>0$, not necessarily of product type, namely $\chi$ satisfies the assumption: 

\begin{description} 
\item{$(\chi)$} there exists $T>0$ such that $\chi ({\tt t})=0$ for ${\tt t}\not\in [-T,T]^N$ and
$\sum_{{\tt k}\in\Z^N} \chi({\tt u-k})=1$, for every ${\tt u}\in\R^N$.
\end{description}
We point out that, since $\chi$ has compact support, obviously $\chi \in L^1(\R^N)$ (and therefore ($\chi_1$) holds) and satisfies the condition ($\chi_2$) since $\chi$ is bounded and the series reduces to a finite sum.

\begin{prop} \label{conv_p}
Let $f\in \Lambda^{p}\cap R_{loc}(\R^N)$, $1\le p<+\infty$. Then, for every $j=1,\ldots,N$,
$$
\lim_{w\to +\infty }\Vert K_{w,j} f-f\Vert_p =0.
$$
\end{prop}

\noindent {\bf Proof.} By assumption $(\chi)$ there holds, for $w>0$, $j=1,\ldots,N$ and ${\tt t}\in\R^N$,
\begin{align*}
| (K_{w,j} f)({\tt t})-f({\tt t})| &= \left| \sum_{{\tt k}\in\Z^N} w\int_{k_j\over w}^{k_{j}+1\over w}f\left({ k_1\over w},\ldots, u\ldots, {k_N\over w}\right) \, du \ \chi(w{\tt t}-{\tt k}) - f({\tt t})\right| \\ 
&= \left| \sum_{{\tt k}\in\Z^N} \left\{w\int_{k_j\over w}^{k_{j}+1\over w}f\left({ k_1\over w},\ldots, u\ldots, {k_N\over w}\right)\,du -f({\tt t}) \right\} \chi(w{\tt t}-{\tt k})\right| 
\\
&\le \sum_{{\tt k}\in\Z^N} \left|w\int_{k_j\over w}^{k_{j}+1\over w}f\left({ k_1\over w},\ldots, u\ldots, {k_N\over w}\right)\,du -f({\tt t}) \right| | \chi(w{\tt t}-{\tt k})| .
\end{align*}
Since $\chi$ has support contained in $[-T,T]^N$, $\chi(w{\tt t}-{\tt k})=0$ if $|w{\tt t}-{\tt k}|>T$, and therefore the series reduces to a finite sum over the indexes ${\tt k}\in\Z^N$ such that 
$(w{\tt t}-{\tt k})\in [-T,T]^N$, namely $\left|{t}_i-{{k}_i \over w}\right| \le {T \over w}$, for every $i=1,\ldots,N$: for such ${\tt k}$
we have that 
$$
\left|w\int_{k_j\over w}^{k_{j}+1\over w}f\left({ k_1\over w},\ldots, u\ldots, {k_N\over w}\right)\,du -f({\tt t}) \right| 
\le \omega_1\left(f;{\tt t}, {2T\over w}, M(\R^N)\right),
$$
by the definition of the modulus of smoothness. Therefore 
\begin{align*}
| (K_{w,j} f)({\tt t})-f({\tt t})| & \le \sum_{(w{\tt t}-{\tt k})\in[-T,T]^N} \left|w\int_{k_j\over w}^{k_{j}+1\over w}f\left({ k_1\over w},\ldots, u,\ldots, {k_N\over w}\right)\,du -f({\tt t}) \right| |\chi(w{\tt t}-{\tt k})| 
\\ &\le \omega_1\left(f;{\tt t}, {2T\over w}, M(\R^N)\right) \sum_{(w{\tt t}-{\tt k})\in[-T,T]^N} |\chi (w{\tt t}-{\tt k})|  \\ &\le A_{\chi} \omega_1\left(f;{\tt t},  {2T\over w}, M(\R^N)\right),
\end{align*}
by ($\chi_2$). Passing to the $L^p-$norm we obtain 
$$
\Vert K_{w,j} f-f \Vert_p \le A_{\chi} \tau_1\left(f; {2T\over w}, M(\R^N)\right)_p,
$$
and the thesis follows by Proposition \ref{tau_p_conv}. \hfill\qed

\section{Convergence in variation for the multidimensional generalized sampling series} \label{sec5}

We are now ready to prove the main result of the paper, that is, the convergence in variation for the multidimensional generalized sampling series with product kernels of averaged type. Since we will use results of the previous section, we assume here that $\chi$ is a kernel with compact support, that is, $\chi:\R^N\rightarrow \R$ satisfies assumption $(\chi)$.

\begin{thm} \label{conv_var}
Let $f\in AC(\R^N)$ be such that ${\partial f\over \partial x_j} \in \Lambda_1 \cap R_{loc} (\R^N)$, for every $j=1,\ldots,N$. 
Then, for every $m\in\N$,
$$
\lim_{w\to +\infty} V[\bar S^m_w f-f]=0.
$$
\end{thm}

\noindent {\bf Proof.} Since $f\in AC(\R^N)$, by Proposition \ref{var-dim},
$\bar S^m_w f\in AC(\R^N)$, for every $w>0$ and $m\in\N$: then 
$$
V[\bar S^m_w f-f] =\int_{\R^N} |\nabla \bar S^m_w f({\tt t}) -\nabla f({\tt t})| \,d{\tt t}.
$$
Now, by Proposition \ref{prop-der}, 
\begin{align*}
V[\bar S^m_w f-f] & \le \int_{\R^N} \sum_{j=1}^N \left| {\partial \over \partial t_j} (\bar S^m_w f)({\tt t}) -{\partial f\over \partial t_j} ({\tt t})\right| \,d{\tt t}\\ 
&=  \sum_{j=1}^N  \int_{\R^N} \left| {1\over m}\sum_{i=1}^m \left(K_{j,w} {\partial f \over \partial t_j}\right)\left({\tt t}'_j,t_j-{m-2(i-1)\over 2w}\right) -{\partial f\over \partial t_j} ({\tt t})\right| \,d{\tt t} \\
&\le  \sum_{j=1}^N {1\over m} \int_{\R^N} \left| \sum_{i=1}^m \left(K_{j,w} {\partial f \over \partial t_j}\right)\left({\tt t}'_j,t_j-{m-2(i-1)\over 2w}\right) \right. \\ &- \left. \sum_{i=1}^m {\partial f \over \partial t_j}\left({\tt t}'_j,t_j-{m-2(i-1)\over 2w}\right) \right| \,d{\tt t}  \\ &+ \sum_{j=1}^N {1\over m} \int_{\R^N} \left| \sum_{i=1}^m {\partial f \over \partial t_j}\left({\tt t}'_j,t_j-{m-2(i-1)\over 2w}\right)  -{\partial f\over \partial t_j} ({\tt t})\right| \,d{\tt t} \\
 &\le  \sum_{j=1}^N  \left\Vert K_{w,j} {\partial f \over \partial t_j} - {\partial f \over \partial t_j} \right\Vert_1 +\\ &+ 
{1\over m} \sum_{j=1}^N \sum_{i=1}^m \int_{\R^N} \left|  {\partial f \over \partial t_j}\left({\tt t}'_j,t_j-{m-2(i-1)\over 2w}\right)-{\partial f\over \partial t_j} ({\tt t})\right| \,d{\tt t}  \\
&:= \sum_{j=1}^N  L_j+ 
{1\over m} \sum_{j=1}^N \sum_{i=1}^m I^i_j.
\end{align*}
Now, $L_j \rightarrow 0$ as $w\rightarrow +\infty$ by Proposition \ref{conv_p}, while $I^i_j\rightarrow 0$ as $w\rightarrow +\infty$, for every $i=1,\ldots,m$, $j=1,\ldots,N$, by the continuity in $L^1$ of the translation operator. Therefore the theorem is proved. \hfill\qed


\section{Examples and applications to smoothing for digital image processing}  \label{sec6}

In this section, we provide some basic examples of kernels for which the above results hold, and we discuss some applications of the variation diminishing-type property (\ref{variation-diminishing}) to smoothing in digital image processing.

  As stated in Section \ref{sec-not}, in this paper we consider product kernels $\bar \chi_m$ of averaged type.
  
\begin{figure}
\centering
\includegraphics[scale=0.28]{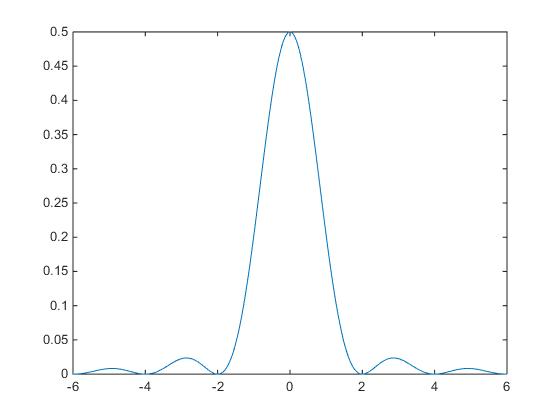}
\hskip0.4cm
\includegraphics[scale=0.28]{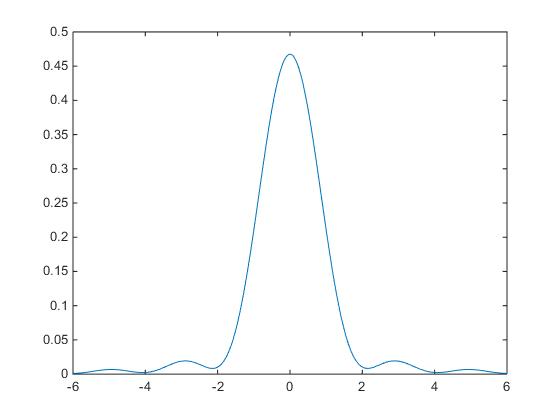}
\caption{On the left: the Fej\'er kernel $F$. On the right: the averaged Fej\'er kernel with $m=1$.} \label{fig1}
\end{figure}

\begin{figure}
\centering
\includegraphics[scale=0.5]{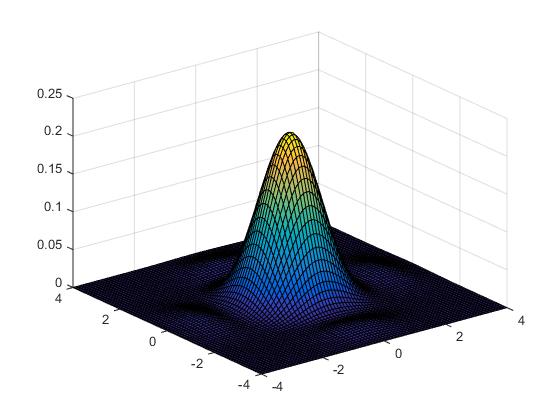}
\caption{The bivariate averaged Fej\'er kernel with $m=1$.} \label{fig2}
\end{figure}

  As a first example, we can consider the multivariate product kernel of the averaged type generated by the well-known Fej\'er kernel (see, e.g., Fig. \ref{fig1} left and \cite{COMIVI1}), defined by:
$$
F(x)\ :=\ {1 \over 2}\, \mbox{sinc}^2(x/2), \hskip1cm x \in \R,
$$
where the sinc-function (see, e.g., \cite{COGA1,ANCOVI2}) is of the form:
$$
\mbox{sinc}(x) := \left\{
\begin{array}{l}
\sin(\pi x) / \pi x, \hskip0.5cm x \neq 0,\\
1,\ \hskip2.1cm x=0.
\end{array}
\right.
$$
Now, the one-dimensional averaged Fej\'er kernel (see, e.g., Fig. \ref{fig1}, right) is the following:
$$
\bar F_m(t)\ :=\ \frac{1}{2\, m} \int_{-m/2}^{m/2}\mbox{sinc}^2\left(\frac{t+v}{2}\right)\, dv, \quad \quad t \in \R, \quad m \in \N,
$$
and the corresponding multivariate version (see, e.g., Fig. \ref{fig2} for the case $N=2$) is:
$$
{\mathcal F}_m({\tt t})\ :=\ \prod_{i=1}^N \bar F_m(t_i), \quad \quad {\tt t} \in \R^N.
$$  
It is well-known that the Fej\'er kernel has unbounded support and satisfies assumptions $(\chi_1)$ and $(\chi_2)$ (see e.g., \cite{COVI4,ANCOVI3}), then for the multivariate generalized sampling series based upon ${\mathcal F}_m$ holds the variation diminishing type property established in Proposition \ref{var-dim}, when multivariate signals of bounded variation are considered.

Moreover, also the relation established in Proposition \ref{prop-der} holds, where the $j$-th first order partial derivatives of the generalized sampling series of a given absolutely continuous function $f$ is related with $K_{w,j} \frac{\partial f}{\partial t_j}$, i.e., the corresponding sampling series of the Kantorovich type of the $j$-th first order partial derivatives of $f$.

  However, the convergence results proved in Section \ref{sec4} can not be applied to the generalized sampling series based upon ${\mathcal F}_m$ since the kernels ${\mathcal F}_m$ do not have compact support, hence assumption $(\chi)$ is not satisfied. 
  
  Other examples of one-dimensional kernels with unbounded support that can be used to define product averaged type kernels can be found, e.g., in \cite{ORTA1,COVI1,COVI2,COVI3}.
  
\begin{figure}
\centering
\includegraphics[scale=0.28]{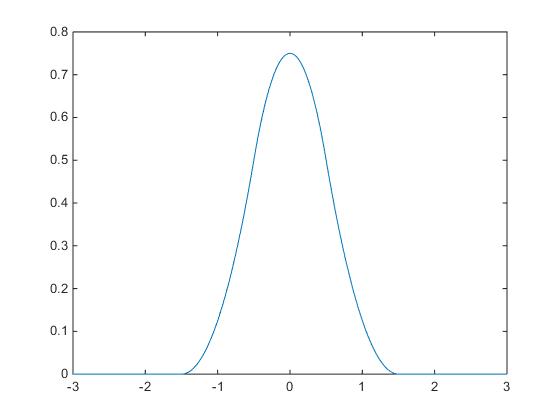}
\hskip0.4cm
\includegraphics[scale=0.28]{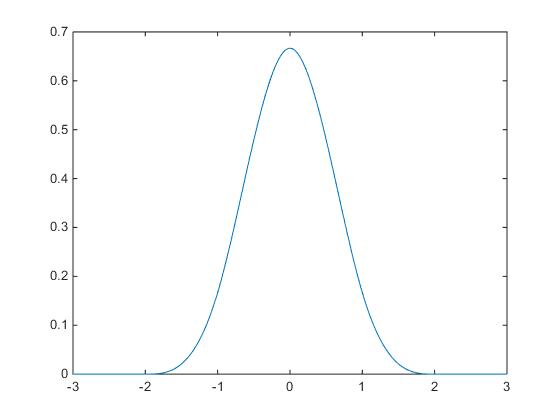}
\caption{On the left: the central B-spline $M_2$. On the right: the central B-spline $M_3$ which corresponds to $\bar M_{2,1}$.} \label{fig3}
\end{figure}

  In order to recall examples of kernels of one variable such that the corresponding multivariate averaged type versions also satisfy assumption $(\chi)$, we recall the definition of the well-known central B-spline of order $n \in \N$ (see, e.g., Fig. \ref{fig3} and \cite{UN1,CHCO1,ALCADE1,ALCADE2,BLRS1}), defined by:
\begin{equation} 
 M_n(x)\ :=\ \frac{1}{(n-1)!} \sum^n_{i=0}(-1)^i \binom{n}{i} 
       \left(\frac{n}{2} + x - i \right)^{n-1}_+,     \hskip0.5cm   x \in \R,
\end{equation}
where $(x)_+ := \max\left\{x,0 \right\}$ denotes ``the positive part'' of $x \in \R$ (see e.g., \cite{ROSA1,DA1}).
The functions $M_n(x)$ are non-negative, continuous with compact support contained in $[-n/2,n/2]$, and satisfy conditions $(\chi_1)$ and $(\chi_2)$.

Now, let us denote by
\begin{equation*}
\bar M_{n,m}(t)\ :=\ m^{-1}\int_{-m/2}^{m/2}M_n(t+v)\, dv, \ \ t \in \R,
\end{equation*}
the averaged B-spline kernel of order $n \in \N$. Recalling the following well-known property:
$$
M_n'(t)\ =\ M_{n-1}(t+1/2)\, -\, M_{n-1}(t-1/2), \ \ t \in \R, \ \ (n \geq 2)
$$
for $m=1$, we have: 
$$
\bar M_{n,1}'(t)\ =\ M_{n}(t+1/2)\, -\, M_{n}(t-1/2)\ =\ M_{n+1}'(t), \ \ t \in \R, \ \  (n \geq 1),
$$
i.e., $\bar M_{n,1}(t)=M_{n+1}(t)+k$, $k \in \R$. Now, since $\bar M_{n,1}$ belongs to $L^1(\R)$ (see \cite{ACV1}), we must have $k=0$ and therefore we conclude that
\begin{equation*}
\bar M_{n,1}(t)\ =\ M_{n+1}(t), \ \  t \in \R, 
\end{equation*}
for every $n \in \N$, namely, the averaged kernel with $m=1$ generated by a central B-spline of order $n$ is a B-spline itself of order $n+1$. 

\begin{figure}
\centering
\includegraphics[scale=0.5]{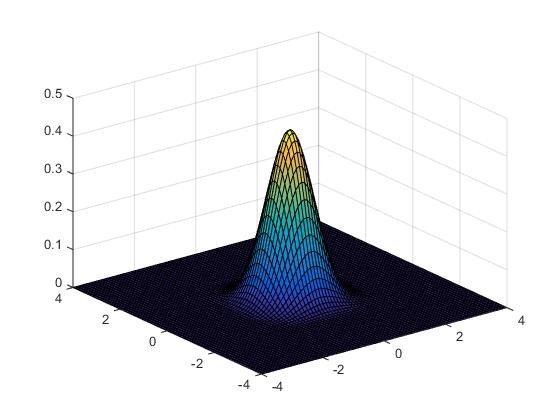}
\caption{The bivariate averaged B-spline kernel of order $2$ with $m=1$ which coincides with ${\cal M}^2_1({\tt t}) = \prod^2_{i=1}M_{3}(t_i)$.} \label{fig4}
\end{figure}

In view of the above remark, we can explicitly state that the multivariate averaged type kernel with $m=1$ and generated by $M_n$ (see, e.g., Fig. \ref{fig4} for the case $n=2$ in two space dimension) is the following:
$$
{\cal M}^n_1({\tt t})\ :=\ \prod^N_{i=1}\bar M_{n,1}(t_i)\ =\ \prod_{i=1}^N M_{n+1}(t_i), \quad {\tt t} \in \R^N.
$$
In practice, in the latter case the multivariate generalized sampling operators based upon averaged B-spline $M_n$ with $m=1$ coincide with the usual generalized sampling series based upon the multivariate central B-spline of order $n+1$.

Generally speaking, we can also define the multivariate average central B-spline kernels as follows:
$$
{\cal M}^n_m({\tt t})\ :=\ \prod^N_{i=1}\bar M_{n,m}(t_i)\ =\ m^{-1}\, \prod_{i=1}^N \int_{-m/2}^{m/2}M_n(t_i+v)\, dv, \quad {\tt t} \in \R^N.
$$
Since ${\cal M}^n_m$ has compact support it satisfies also assumption $(\chi)$; thus the above sampling series fulfills the results of both Section \ref{sec-mres} and Section \ref{sec5}. Further, for the latter examples of kernels also Proposition \ref{conv_p} and Theorem \ref{conv_var} can be applied.

\vskip0.2cm

Now, at the end of this section we consider some applications to smoothing in digital processing. For basic facts concerning this numerical tool for imaging, see e.g., \cite{BI1,PA1,SO1}.

   It is well-known that any static gray scale image is a bivariate signal with compact support; then it can be (naturally) modeled as follows:
$$
I_A(x,y)\ :=\ \sum^{m}_{i=1}\sum^{m}_{j=1}a_{ij} \cdot \textbf{1}_{ij}(x,y) \hskip1cm ((x,y) \in \R^2),
$$
for every image (matrix) $A=(a_{ij})_{ij}$, $i,j=1,2,...,m$, where $\textbf{1}_{ij}(x,y)$, $i,j =1,2,...,m$, is the characteristic function of the sets $(i-1,\ i]\times(j-1,\ j]$ (i.e. $\textbf{1}_{ij}(x,y)=1$, for $(x,y)\ \in\ (i-1,\ i]\times(j-1,\ j]$ and $\textbf{1}_{ij}(x,y)=0$ otherwise).

Note that the above function $I_A(x,y)$ is defined in such a way that to every pixel $(i,j)$ it is associated the corresponding gray level $a_{ij}$.

 Moreover, by the above representation $I_A$ of the image $A$ it turns out that $I_A \in BV(\R^2)$, hence one can consider approximations of $A$ by means of the bivariate generalized sampling series based upon the averaged type kernels. The main advantage that can be achieved by the above procedure is expressed by Proposition \ref{var-dim}: $\bar S^m_w I_A$ have no-bigger variation than $I_A$, for sufficiently large $m \in \N$, i.e., the operators provide an approximated version of the original image $A$, but regularized. In this sense, the generalized sampling series can be used for smoothing of images. 
 
   Clearly, in order to visualize an approximation (new image) of the original image $A$ by means of $\bar S^m_w I_A$, we need to sample the operators, for  $w>0$, with a fixed sampling rate concordant with the dimension of $A$. Obviously, the sampling rate is chosen arbitrarily hence one can also consider different (high) sampling rates.

    The effect of the proposed procedure can be strongly noticed at the edges of the figures, where the jumps of gray levels are reduced with respect to the corresponding ones in the original image.

   Now we can give the following practical examples of image reconstruction in order to show the smoothing capabilities of the above operators. An optimized version of the above described algorithm for image reconstruction and smoothing can be implemented by means of the MATLAB programming language, following the indications outlined in \cite{ING1,ING2} in case of the so-called sampling Kantorovich algorithm for image enhancement.   
   
\begin{figure}
\centering
\includegraphics[scale=3]{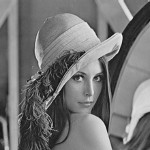}
\hskip0.75cm
\includegraphics[scale=1.2]{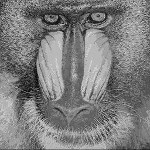}
\caption{On the left: Lena. On the right: Baboon.} \label{fig5}
\end{figure}

   For the numerical experiments, we consider the well-known images of Lena and Baboon with $150\times 150$ pixel resolution (see Fig. \ref{fig5}). 
\begin{figure}
\centering
\includegraphics[scale=0.7]{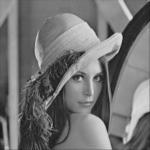}
\hskip0.75cm
\includegraphics[scale=0.7]{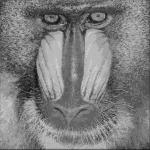}
\caption{On the left: the reconstruction of Lena (of $150\times 150$ pixel) by means of the operator $\bar S^m_w$, $w=4$, based upon the bivariate averaged Fej\'er kernel with $m=4$. On the right: the reconstruction of Baboon made as for Lena.} \label{fig6}
\end{figure}

   First of all, we reconstruct the original images in Fig. \ref{fig5} by using the bivariate averaged Fej\'er kernel with $m=4$. According to Proposition \ref{var-dim}, and observing that $\|F\|_1=1$, we have:
\begin{equation}  \label{disug-variaz}
V[\bar S^4_w f]\ \le\ {1\over 4}\, V[f],
\end{equation}
i.e., the total variation of any reconstructed images is, at least, 4 times smaller than that of the original ones, producing the smoothing effect. In Fig. \ref{fig6} we have the reconstruction of Lena and Baboon (of $150\times 150$ pixel resolution) by means of the operator $\bar S^m_w$, $w=4$, based upon the bivariate averaged Fej\'er kernel with $m=4$. By detailed analysis of the edges (especially in case of Lena, at the contours of the hat) it is possible to observe the smoothing of the analyzed images.
\begin{figure}
\centering
\includegraphics[scale=0.7]{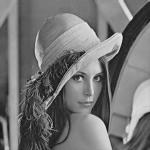}
\hskip0.75cm
\includegraphics[scale=0.7]{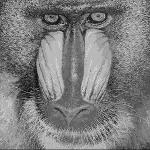}
\caption{On the left: the reconstruction of Lena (of $150\times 150$ pixel) by means of the operator $\bar S^m_w$, $w=10$, based upon the bivariate averaged Fej\'er kernel with $m=4$. On the right: the reconstruction of Baboon made as for Lena.} \label{fig7}
\end{figure}
In Fig. \ref{fig7} we have the reconstruction of Lena and Baboon (of $150\times 150$ pixel resolution) by means of the operator $\bar S^m_w$, $w=10$, based upon the bivariate averaged Fej\'er kernel with $m=4$.

  The main differences that can be observed between the images in Fig. \ref{fig6} and Fig. \ref{fig7} are that, for big values of $w$ the images are closer to the original and consequently the edges tend to be more clear. 
\begin{figure}
\centering
\includegraphics[scale=0.7]{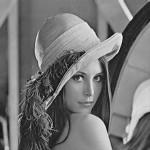}
\hskip0.75cm
\includegraphics[scale=0.7]{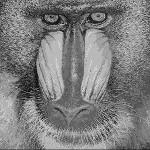}
\caption{On the left: the reconstruction of Lena (of $150\times 150$ pixel) by means of the operator $\bar S^m_w$, $w=6$, based upon the bivariate averaged central B-spline of order $3$ with $m=4$. On the right: the reconstruction of Baboon made as for Lena.} \label{fig8}
\end{figure}
Finally, in Fig. \ref{fig8} we have the reconstruction of Lena and Baboon (of $150\times 150$ pixel resolution) by means of the operator $\bar S^m_w$, $w=6$, based upon the bivariate central B-spline of order $3$ with $m=4$. Note that, also for the latter case one can state a inequality as that given in (\ref{disug-variaz}) since also $\|M_3\|_1=1$.


\section*{Acknowledgments}

The authors are members of the Gruppo  
Nazionale per l'Analisi Matematica, la Probabilit\'a e le loro  
Applicazioni (GNAMPA) of the Istituto Nazionale di Alta Matematica (INdAM). 

\noindent The authors are partially supported by the "Department of Mathematics and Computer Science" of the University of Perugia (Italy). The first and the third author have been partially supported within the project "Metodi di teoria degli operatori e di Analisi Reale per problemi di approssimazione ed applicazioni", funded by the 2017 basic research fund of the University of
Perugia. Finally, the second author of the paper has been partially supported within a 2018 GNAMPA-INdAM Project: ``Dinamiche non autonome, analisi reale e applicazioni''.

\end{document}